\documentclass[11pt]{article}
\usepackage{amssymb,latexsym,amsfonts,verbatim,amscd}
\usepackage[all]{xy}
\newtheorem{Def}{Definition}[section]
\newtheorem{Thm}[Def]{Theorem}

\newtheorem{Prop}[Def]{Proposition}
\newtheorem{Cor}[Def]{Corollary}
\newtheorem{Rem}[Def]{Remark}
\newtheorem{Que}[Def]{Question}

\newtheorem{Fac}[Def]{Fact}
\newtheorem{Ex}[Def]{Example}

\font\nat msbm10 scaled\magstephalf
\def\N{\hbox{\nat\char78}}

\def\Q{\hbox{\nat\char81}}

\def\telos{\hfill$\dashv$}
 \font\goth eufm10

\def\int{\hbox{\rm int}}

\begin{document}
\sloppy

\title{Russell's typicality as another randomness notion}
\author{Athanassios Tzouvaras}

\date{}
\maketitle

\begin{center}
Aristotle University of Thessaloniki \\
Department  of Mathematics  \\
541 24 Thessaloniki, Greece. \\
e-mail: \verb"tzouvara@math.auth.gr"
\end{center}

\begin{abstract}
We reformulate slightly  Russell's notion of typicality, so as to  eliminate its circularity and make it applicable to elements of any first-order structure. We argue that the notion parallels Martin-L\"{o}f (ML) randomness, in the sense that it uses definable sets  in place of computable ones and sets of ``small'' cardinality (i.e., strictly smaller than that of the structure domain) in place of measure zero sets. It is shown that if the domain $M$ satisfies $cf(|M|)>\aleph_0$, then there exist $|M|$ typical elements and only  $<|M|$ non-typical ones. In particular this is true for the standard model ${\cal R}$ of second-order arithmetic.
By allowing parameters in the defining formulas, we are led to relative typicality, which satisfies most of van Lambalgen's  axioms for relative randomness.  However van Lambalgen's theorem is false for relative typicality. The class of typical reals is incomparable (with respect to $\subseteq$) with the classes of  ML-random, Schnorr random and computably  random  reals. Also the class of  typical reals is closed under Turing degrees and under the jump operation (both ways).
\end{abstract}

\vskip 0.2in


\vskip 0.2in


\section{Introduction}
Bertrand Russell \cite[p. 89]{Ru95}, in an attempt to explain  impredicative definitions (and the need to appeal to the reducibility axiom  in order to avoid them), puts under examination the following definition  of  ``typical Englishman'':  {\em ``Suppose, for example, you were to  suggest that `a typical Englishman is one who possesses} all {\em the properties possessed by a majority of Englishmen'.''} Now typicality so defined is itself a property of
Englishmen, so  in order for someone to check whether a particular Englishman
is typical, one has to check, among other things, if  he  is
typical, so checking enters a circle and fails.\footnote{Russell  seems to worry not about circularity itself but rather about a  contradiction which  is supposed to emerge from this. For he continues the above phrase as follows: {\em ``You will easily realise that most Englishmen do not possess}  all {\em   the properties that most Englishmen possess, and therefore a typical Englishman, according to your own definition, would be untypical.''} (ibid.) It is not clear to me how Russell concludes that {\em  most Englishmen do not possess  all the properties that most Englishmen possess,} and thus how the contradiction is derived. I suppose that Russell's conclusion {\em ``You will easily realise that...''} is empirical rather than logical. In section 2 below we show (Example \ref{Ex:rationals} and  Theorems \ref{T:typicalgen}, \ref{T:typicalreal}) that there are plenty of   structures in which  the {\em majority} of elements, or even the {\em totality} of them, can be typical. So in these structures,  the property of metalanguage expressing typicality  is  itself typical. }

Despite this failure, the core of the definition remains natural and  sound, and  can be easily restored  if we  make some   changes that will eliminate circularity. It is  the purpose of this paper to provide such a  strict definition of typicality, and then explore its behavior in concrete contexts and its connections with randomness.  We believe  that the intended meaning of  Russell's notion can be described as follows: given a structured universe of things, which could  mathematically be represented by  a first-order structure $\mbox{\goth M}=(M,\ldots)$,  an element  $a$ of $M$ is  typical if it is {\em not special} in any {\em definable} sense, that is,  if it does not possess any  property  that makes it belong to a {\em special definable minority} $X\subset M$. This can be done rigorously in three steps: (a) Define precisely when a set $X\subseteq M$ contains a  majority/minority of elements  of $M$.  (b) Specify what a ``typical property'' (of the first-order language $L$ of $\mbox{\goth M}$) is in terms of the previously specified majority notion.   (c) Define an element $a$ of $M$ to be ``typical'' if it satisfies all typical properties. Since typical properties are properties of the object language $L$, while  typicality  is  a property of the metalanguage, no circularity occurs.   In  section 2 we make these steps  precise.

In section 3 we argue that  typicality is a kind  of randomness notion, especially when we apply it to the standard model of full second-order arithmetic  ${\cal R}=(\omega,{\cal P}(\omega),+,\cdot,<,\in,0,1)$.  In the context  of  reals, the most popular  notion of randomness is Martin-L\"{o}f (ML-) randomness (\cite{DH08} is a rich source of information about this notion and its  variants). The parallelism between ML-random and typical reals lies on this:

$\bullet$ The definition of a ML-random real is based on the  combination of computability and small sets in the sense of measure theory (sets of measure zero). Namely, $a$ is ML-random if it cannot be ``trapped'' within any set $\bigcap_nV_n$ of measure zero, where $(V_n)_n$ is a computable sequence of r.e. open sets such that $\mu(V_n)\leq 2^{-n}$.

$\bullet$ The definition of a  typical real is based on the combination of definability and small  sets in the sense of cardinality (sets of cardinality strictly less than the continuum). Namely, $a$ is typical if it cannot be ```trapped'' within any set $X$ which is definable and small, i.e., $|X|<2^{\aleph_0}$.

Thus typicality arises if we  replace computable sets with definable ones and small sets from the point of view of measure theory with small sets from the point of view of  cardinality.

There are structures, like $(\N,S,+,\cdot,0)$, that contain no typical elements, and others, like $(\Q,<)$, all elements of which are  typical. We give also some sufficient conditions in order for a structure to contain typical elements. One such  simple and basic condition is given in  Theorem \ref{T:typicalgen}. It  says  that  every $\mbox{\goth M}=(M,\ldots)$ with a countable language   such that $cf(|M|)>\aleph_0$ contains $|M|$ typical elements, while  only $<|M|$ non-typical ones.  In particular this is true for the structure ${\cal R}$  of reals mentioned above, since  $cf(2^{\aleph_0})>\aleph_0$.

If further  we allow parameters to occur in properties, we are led to the notion of relative  typicality ``$x$ is typical with respect to $\overline{y}$'', denoted $\textsf{Tp}(x,\overline{y})$, which is analogous to the relative randomness $R(x,\overline{y})$ studied  by van Lambalgen in \cite{La90}, \cite{La92} and \cite{La96}. We show that  most of the  basic axioms for  $R(x,\overline{y})$ considered in these papers hold also for $\textsf{Tp}(x,\overline{y})$ in general structures, and in particular in the structure ${\cal R}$. On the other hand the two notions deviate significantly at certain points. For example van Lambalgen's theorem does not hold for relative typicality, while it holds for ML-randomness. Actually the class of typical reals is incomparable (with respect to inclusion) with the classes of  ML-random, Schnorr random and computably random  reals.

Finally we show that the relation $\textsf{Tp}(x,\overline{y})$  is closed with respect to Turing reducibility $\leq_T$ and Turing degrees, as well as with respect to the jump operator $a\mapsto a'$.

Some open questions are stated at the end of the paper.

\section{Formalizing typicality}
Let $L$ be  a first-order language, $\mbox{\goth M}=(M,\ldots)$  an $L$-structure and $A\subseteq M$. $L(A)$ denotes $L$ augmented with parameters from $A$. By some abuse of language we refer also to $L(A)$ as the set of formulas of $L(A)$.   By a property of $L(A)$ we mean  a formula $\phi(x)\in L(A)$ with  one free variable. We have first to make precise what it means for  a set $X\subseteq M$ to contain the majority of elements of $M$.  The definition is the expected one: a majority subset is any  set that contains strictly more elements than its complement, i.e., any  $X\subseteq M$ such that  $|X|>|M\backslash X|$. We shall refer to such a $X$ as a {\em majority subset} (and accordingly to $M\backslash X$ as a {\em minority subset}). The definition applies both to finite and infinite domains  $M$. If $M$ is finite, then $X$ is a majority set just if  $|X|>|M|/2$.  If  $M$ is infinite (as will be the case throughout this paper), the above condition is equivalent to  $|M\backslash X|<|M|$ (which in particular  implies that $|X|=|M|$). Let us  set
$$\textsf{mj}(M)=\{X\subseteq M:|X|>|M\backslash X|\}=\{X:|M\backslash X|<|M|\}$$
for the class of majority subsets of $M$.
It is easy to see that $\textsf{mj}(M)$ is a filter on $M$, which extends the Fr\'{e}chet filter of cofinite sets
$$Fr(M)=\{X\subseteq M:|M\backslash X| \ \mbox{is finite}\}.$$
For countable $M$,  $\textsf{mj}(M)=Fr(M)$ and  $\textsf{mj}(M)$ is sometimes  referred to as {\em generalized Fr\'{e}chet} filter.

Having fixed a rigorous notion of majority, the  formalization of Russell's concept  comes in  two steps:  first, define  a property $\phi(x)$ of a language $L$ for a structure  $\mbox{\goth M}=(M,\ldots)$ to be {\em typical}, if it  defines a majority subset of $M$, i.e., if its extension in $\mbox{\goth M}$ belongs to $\textsf{mj}(M)$; second, define an  element $a\in M$ to be {\em typical}, if it satisfies all typical properties for  $\mbox{\goth M}$. A weak notion of typicality is obtained (over uncountable structures) if we use the filter $Fr(M)$ in  place of $\textsf{mj}(M)$.

\subsection{Typical and weakly typical properties and elements of a  first-order structure. Relative typicality}
Given an $L$-structure $\mbox{\goth M}$, a set $A\subseteq M$  and a property $\phi(x)$ of $L(A)$, we often denote by $ext(\phi)^{\mbox{\goth M}}$, or just $ext(\phi)$, the extension of $\phi(x)$ in $\mbox{\goth M}$, i.e., $$ext(\phi)=\{a\in M:\mbox{\goth M}\models\phi(a)\}.$$

\begin{Def} \label{D:typical}
{\em    A property $\phi(x)$ of $L(A)$ is said to be}  typical over $\mbox{\goth M}$, {\em or} $(\mbox{\goth M}, A)$-typical, {\em or just} $A$-typical {\em (resp.} weakly $A$-typical,{\em )  if $ext(\phi)\in \textsf{mj}(M)$ (resp. $ext(\phi)\in Fr(M)$). In particular  $\phi(x)$ is} typical {\em (resp.} weakly typical{\em ) if it is $\emptyset$-typical (resp. weakly $\emptyset$-typical).}\footnote{An equivalent way to define typicality and weak typicality of properties is by the use of  generalized quantifiers $Q_{most}$ and $Q_{in\!f}$, respectively.  For every structure $\mbox{\goth M}$, the interpretation of  $Q_{most}$  and $Q_{in\!f}$ is  the majority filter $\textsf{mj}(M)$ and the Fr\'{e}chet filter, respectively. That is, for every $\phi(x)$ of $L(A)$,
$$\mbox{\goth M}\models (Q_{most}x)\phi(x) \ \iff \ \{a\in M: \mbox{\goth M}\models \phi(a)\}\in \textsf{mj}(M),$$
and
$$\mbox{\goth M}\models (Q_{in\!f}x)\phi(x) \ \iff \ \{a\in M: \mbox{\goth M}\models \phi(a)\}\in Fr(M).$$
Then a property $\phi(x)$ is $A$-typical (resp. weakly $A$-typical) over $\mbox{\goth M}$ iff $\mbox{\goth M}\models (Q_{most}x)\phi(x)$ (resp. $\mbox{\goth M}\models (Q_{in\!f}x)\phi(x)$).}
\end{Def}
Let $TP(\mbox{\goth M},A)$, $wTP(\mbox{\goth M},A)$ denote the sets   of  $A$-typical and weakly $A$-typical   properties $\phi(x)$ of $L(A)$ over $\mbox{\goth M}$. Obviously $wTP(\mbox{\goth M},A)\subseteq TP(\mbox{\goth M},A)$.

\begin{Fac} \label{F:type}
Both $TP(\mbox{\goth M},A)$ and $wTP(\mbox{\goth M},A)$ are types  over $\mbox{\goth M}$ with parameters in $A$ in the usual sense, i.e., they are  finitely satisfiable in $\mbox{\goth M}$.
\end{Fac}

{\em Proof.} It is straightforward from the definition that the conjunction of finitely many $A$-typical (resp. weakly $A$-typical) properties is $A$-typical (resp. weakly $A$-typical), thus satisfiable.   \telos

\vskip 0.2in

We shall be mostly interested in $A$-typical properties and elements  for finite $A$. In this case the elements of $A$ occur in $\phi$ in a certain order, i.e., as a vector $\overline{a}=(a_1,\ldots,a_n)$, so it convenient  to use in parallel the vector notation and write  also $\overline{a}$-typical (weakly $\overline{a}$-typical) instead of $A$-typical (weakly $A$-typical).
Concerning the existence of $A$-typical properties over a structure $\mbox{\goth M}$, we have the following simple facts. All structures $\mbox{\goth M}$ considered below are infinite. The proof of the following is straightforward.

\begin{Fac} \label{F:proeprtyexistence}
(i) Every tautology $\phi(x)$ (e.g. $\psi(x)\vee\neg\psi(x)$) is $M$-typical.

(ii) For  every  tuple $\overline{a}=(a_1,\ldots,a_n)$ of $M$, the property
$$\phi_{\overline{a}}(x):=(x\neq a_1)\wedge\cdots\wedge (x\neq a_n)$$
is $\overline{a}$-typical.

(iii) If  $\overline{b}$ is a tuple of $\overline{a}$-definable elements of $M$, then the property $\phi_{\overline{b}}(x)$ is equivalent over $\mbox{\goth M}$ to an $\overline{a}$-typical property.
\end{Fac}

We come to typical elements.

\begin{Def} \label{D:typicalel}
{\em Given a structure $\mbox{\goth M}$ and $A\subseteq M$,  an element $a\in M$ is said to be} $(\mbox{\goth M}, A)$-typical, {\em or just} $A$-typical {\em (resp.} weakly $A$-typical{\em ), if it satisfies the type  $TP(\mbox{\goth M},A)$ of $A$-typical properties over $\mbox{\goth M}$ (resp. the type $wTP(\mbox{\goth M},A)$ of weakly $A$-typical properties), i.e $\mbox{\goth M}\models \phi(a)$ for every $\phi\in TP(\mbox{\goth M},A)$. }
\end{Def}

Allowing parameters to be used in definitions,  gives rise to the following {\em relative typicality} relations  between elements of a structure $\mbox{\goth M}$:

$\bullet$ $\textsf{Tp}(a,\overline{b})$: ``$a$ is $\overline{b}$-typical'',

$\bullet$ $\textsf{wTp}(a,\overline{b})$:  ``$a$ is weakly $\overline{b}$-typical''.

In particular we write $\textsf{Tp}(a)$,  $\textsf{wTp}(a)$  instead of $\textsf{Tp}(a,\emptyset)$, $\textsf{wTp}(a,\emptyset)$, respectively.

Recall that given an $L$-structure $\mbox{\goth M}$ and $A\subseteq M$, an element $b\in M$ is said to be {\em $A$-definable} in $\mbox{\goth M}$, if there is a formula $\phi(x,\overline{y})$ of $L$ and a tuple $\overline{a}\in A$ such that $b$ is the unique element of $M$ that satisfies $\phi(x,\overline{a})$, i.e.,
$$\mbox{\goth M}\models(\forall x)(x=b\leftrightarrow \phi(x,\overline{a})).$$
An element is said to be definable if it is $\emptyset$-definable.
More generally $b\in M$ is said to be {\em $A$-algebraic} in $\mbox{\goth M}$, if there is a formula $\phi(x,\overline{y})$ and a tuple $\overline{a}\in A$ such that
$\mbox{\goth M}\models\phi(b,\overline{a})$ and $ext(\phi(x,\overline{a}))$ is finite. For every tuple $\overline{a}\in M$, let
${\rm df}(\overline{a})$ and  ${\rm cl}(\overline{a})$ denote the sets of $\overline{a}$-definable and $\overline{a}$-algebraic elements of $M$, respectively. Obviously ${\rm df}(\overline{a})\subseteq {\rm cl}(\overline{a})$.

The following are immediate consequences of the definitions.

\begin{Fac} \label{F:easy}
For any structure  $\mbox{\goth M}$ and elements of $M$ the following hold:

(i) $\neg \textsf{Tp}(a,\overline{b})$ iff  there exists $\phi(x,\overline{y})$ such that $\mbox{\goth M}\models\phi(a,\overline{b})$ and $|ext(\phi(x,\overline{b}))|<|M|$.

(ii) $\neg \textsf{wTp}(a,\overline{b})$ iff there exists $\phi(x,\overline{y})$ such that $\mbox{\goth M}\models\phi(a,\overline{b})$ and $ext(\phi(x,\overline{b}))$  is finite.

(iii) If $\{b_1,\ldots b_n\}\subseteq \{c_1,\ldots,c_m\}$, then $\textsf{Tp}(a,\overline{c})$ implies $\textsf{Tp}(a,\overline{b})$ and $\textsf{wTp}(a,\overline{c})$ implies $\textsf{wTp}(a,\overline{b})$.

(iv) $\textsf{Tp}(a,\overline{b}) \ \Rightarrow \ \textsf{wTp}(a,\overline{b}) \ \Leftrightarrow \ a\notin {\rm cl}(\overline{b})$

(v) For countable  $\mbox{\goth M}$,
$\textsf{Tp}(a,\overline{b}) \ \Leftrightarrow \ \textsf{wTp}(a,\overline{b})$.
\end{Fac}

While for all (infinite) structures there exist typical properties, one cannot expect that all structures contain typical elements. E.g. this is the case with structures all elements  of which are ``special'', i.e., definable  without parameters. On the other hand,  there exist structures consisting  entirely of typical  elements. Two such prominent  examples of structures, lying at the two opposite ends of the spectrum,  are considered below.

\begin{Ex} \label{Ex:integers}
{\em Let $\mbox{\goth M}={\cal N}=(\N,S,+,\cdot,0)$ be the standard structure of natural numbers.  Every $n\in\N$ is definable in ${\cal N}$, so the formulas of $L(\N)$ coincide with those of $L$. Thus for any tuple $\overline{n}=(n_1,\ldots,n_k)$ of elements of  $\N$, the formula $\phi_{\overline{n}}(x)$ defined in \ref{F:proeprtyexistence} is typical. Coming to elements, it is easy to see that ${\cal N}$ contains no typical element. Because  for every  tuple $\overline{n}=(n_1,\ldots,n_k)$ of $\N$, the typical  property  $\phi_{\overline{n}}(x):=(x\neq n_1)\wedge\cdots\wedge (x\neq n_k)$ has extension  $\N\backslash \{n_1,\ldots,n_k\}$.  So if $a\in \N$ were typical, it should belong to the extensions of all properties $\phi_{\overline{n}}(x)$, i.e. to $\bigcap Fr(\N)$, but this is empty.}
\end{Ex}

\begin{Ex} \label{Ex:rationals}
{\em Let  $\mbox{\goth M}={\cal Q}=(\Q,<)$ be the ordered set of rationals. Here, in full contrast to ${\cal N}$, the elements
of ${\cal Q}$ are all of  the {\em same type,}   i.e., for every $\phi(x)$ and any $a,b\in \Q$, ${\cal Q}\models\phi(a)\leftrightarrow\phi(b)$. So for every $\phi(x)$, either it holds of all elements, i.e., $ext(\phi)=\Q$, so $\phi(x)$ is typical, or it holds of no element, in which case $ext(\neg \phi)=\Q$, i.e., $\neg\phi(x)$ is typical. It follows that for every $\phi(x)$, either $\phi(x)$ or $\neg\phi(x)$ is typical over ${\cal Q}$.

Let  $A\subseteq \Q$ be a set of parameters. The theory $DLO^*$ (of dense linear order without end-points) admits quantifier elimination, so each property  $\phi(x)$ with parameters from $A$ is equivalent over ${\cal Q}$ to a Boolean combination of  formulas $x<a_i$, $a_j<x$, $x=a_k$, and their negations.  Of them only the properties  $x\neq a_i$ are $A$-typical. So the $A$-typical properties over ${\cal Q}$ are exactly $\phi_{\overline{a}}(x):=(x\neq a_1)\wedge\cdots\wedge (x\neq a_n)$
for $\overline{a}\in A$. It follows that  $a$ is $A$-typical iff it satisfies all such $\phi_{\overline{a}}(x)$, i.e., iff it belongs to  $\bigcap_{a_1,\ldots,a_n\in A} \Q\backslash\{a_1,\ldots,a_n\}$ which equals $\Q\backslash A$. This is exactly the set of $A$-typical elements of ${\cal Q}$.}
\end{Ex}

Here is a sufficient condition for the existence of typical elements in a structure $\mbox{\goth M}$.

\begin{Thm} \label{T:suffices}
Let $\mbox{\goth M}$ be a structure. If $\mbox{\goth M}$ is $\kappa$-saturated, for some $\kappa\geq \aleph_0$, then for every $A\subseteq M$ such that $|A|<\kappa$, $M$ contains $A$-typical elements.
\end{Thm}

{\em Proof.} We saw that the set  $TP(\mbox{\goth M},A)$ of $A$-typical properties is a type  over $(\mbox{\goth M},A)$. Sine $|A|<\kappa$ and  $\mbox{\goth M}$ is $\kappa$-saturated, $TP(\mbox{\goth M},A)$ is satisfiable in $\mbox{\goth M}$. Every $a\in M$ satisfying $TP(\mbox{\goth M},A)$ is $A$-typical. \telos

\vskip 0.2in

Here is another sufficient condition for the existence of typical elements,  independent of saturation.

\begin{Thm} \label{T:typicalgen}
Let $\mbox{\goth M}$ be an $L$-structure, for a countable $L$, and let $A\subseteq M$ be a set of parameters  such that $cf(|M|)>\max(\aleph_0,|A|)$.  Then $\mbox{\goth M}$ contains $|M|$ $A$-typical elements, while only $<|M|$ non-$A$-typical ones.
\end{Thm}

{\em Proof.} Let $\mbox{\goth M}$ and $A\subseteq M$ be as stated. By definition an element $a\in M$ is non-$A$-typical iff it satisfies a property with minority extension.  So if $$S=\{\phi(x,\overline{b}):\overline{b}\in A \ \& \ |ext(\phi(x,\overline{b}))|<|M|\},$$
then  the set of all non-$A$-typical elements of $M$ is
$X=\bigcup\{ext(\phi):\phi\in S\}$, and hence
$|X|\leq \Sigma\{|ext(\phi)|:\phi\in S\}$.
Now  $|S|\leq \max(\aleph_0,|A|)<cf(|M|)$, i.e., $|S|<cf(|M|)$, and also   for every $\phi\in S$, $|ext(\phi)|<|M|$. It follows that $\Sigma\{|ext(\phi)|:\phi\in S\}<|M|$, hence $|X|<|M|$. The set of $A$-typical elements is the complement of $X$, $M\backslash X$, so $|M\backslash X|=|M|$.  \telos

\begin{Cor} \label{C:countable}
Let $\mbox{\goth M}$ be an $L$-structure, for a countable $L$, such that $cf(|M|)>\aleph_0$. Then for any  $A\subseteq M$ such that $|A|\leq \aleph_0$,  $\mbox{\goth M}$ contains $|M|$ $A$-typical elements and $<|M|$ non-$A$-typical ones.
\end{Cor}

Let $Z_2$ be the theory of full second-order arithmetic, whose language is $L_2=\{+,\cdot,<,\in,0,1\}$. $L_2$ has  variables $m,n,i,j,\ldots$ for numbers and  variables $x,y,z$  for sets of numbers, i.e., for {\em reals} (for details see \cite{Si99}). The full standard model of $Z_2$ is ${\cal R}=(\omega,{\cal P}(\omega),+,\cdot,<,\in,0,1)$.  $a,b,c,\ldots$ denote arbitrary reals. Since $cf(|{\cal P}(\omega)|)=cf(2^{\aleph_0})>\aleph_0$,  applying Corollary \ref{C:countable} to the structure  ${\cal R}$, we have the following immediate consequence.

\begin{Thm} \label{T:typicalreal}
For every $A\subseteq {\cal P}(\omega)$ such that $|A|\leq \aleph_0$, there exist $2^{\aleph_0}$ $A$-typical reals, while only $<2^{\aleph_0}$ non-$A$-typical ones. More precisely:  For every finite (or even countable) tuple $\overline{b}$ of reals, $|\{x:\textsf{Tp}(x,\overline{b})\}|=2^{\aleph_0}$, while  $|\{x:\neg \textsf{Tp}(x,\overline{b})\}|<2^{\aleph_0}$.
\end{Thm}

Example \ref{Ex:rationals} and Theorems \ref{T:typicalgen}, \ref{T:typicalreal} provide examples of structures  that not only contain typical elements, but are structures over which the property of  typicality itself (as a property of the metalanguage) is also ``typical'' (in the sense that its extension is a majority set).

It follows from  Fact \ref{F:easy} (ii)  that every  algebraic (hence also every definable) element of $\mbox{\goth M}$ is non-typical. However, using  the coding capabilities of  reals, it is easy to see that for any $\overline{b}$, the  $\overline{b}$-algebraic and $\overline{b}$-definable reals coincide.

\begin{Fac} \label{F:coincide}
For all $\overline{b}\in {\cal P}(\omega)$,
${\rm cl}(\overline{b})={\rm df}(\overline{b})$.
\end{Fac}

{\em Proof.} It suffices to see that ${\rm cl}(\overline{b})\subseteq{\rm df}(\overline{b})$. This is a simple consequence of the fact that every finite (and even a countable) set of reals can be coded by a single real. Let $a\in {\rm cl}(\overline{b})$. Then there is a formula $\phi(x,\overline{y})$ of $L_2$ such that $\models\phi(a,\overline{b})$ and $\{x:\models\phi(x,\overline{b})\}$ is finite, so let $\{x:\models\phi(x,\overline{b})\}=\{c_1,\ldots,c_n\}$. A code for the latter set is the real  $c=\{\langle i,m\rangle:m\in c_i, 1\leq i\leq n\}$, where $\langle i,m\rangle$ is the usual arithmetical code of the pair $(i,m)$. If $(c)_i=\{m:\langle i,m\rangle\in c\}$, then  $c_i=(c)_i$ for  $1\leq i\leq n$. In particular $a=(c)_{i_0}$ for some $i_0$. Now $c$ is defined by the formula $\psi(x,\overline{b}):=(\forall y)(\phi(y,\overline{b})\leftrightarrow y=(x)_1\vee\cdots\vee(x)_n)$. Thus  $c\in{\rm df}(\overline{b})$,  and since also $a\in {\rm df}(c)$, it follows  $a\in {\rm df}(\overline{b})$.  \telos

\vskip 0.2in

Fact \ref{F:coincide}, together with \ref{F:easy} (ii), implies:

\begin{Fac} \label{F:reduce}
In the structure ${\cal R}$, for all reals $a,\overline{b}$ we have
$$\textsf{Tp}(a,\overline{b}) \Rightarrow  \textsf{wTp}(a,\overline{b})\Leftrightarrow a\notin {\rm df}(\overline{b}).$$
\end{Fac}

It is well-known that the $L_2$-definable subsets of ${\cal P}(\omega)$ are those of the analytical (or lightface) hierarchy  of descriptive set theory consisting of the classes  $\Sigma_n^1$, $\Pi_n^1$, for $n\geq 0$. So  a real $a$ is non-typical iff it belongs to a minority set of this hierarchy, i.e., iff $a\in X$ for some   $X\in \Sigma_n^1$,  for some $n\geq 0$, such that  $|X|<2^{\aleph_0}$. The problem is considerably simplified if we assume $CH$, in which case $|X|<2^{\aleph_0}$ means  $X$ is countable. Thus we have the following.

\begin{Fac} \label{F:simplify}
Assuming  $CH$, a real $a\in {\cal P}(\omega)$ is non-typical iff there is a countable $X\in \bigcup_n\Sigma^1_n$  such that $a\in X$.
\end{Fac}
Interestingly enough, the {\em countable} sets of the analytical hierarchy have drawn  the attention of important researchers in descriptive set theory quite early, with remarkable results  (see \cite{KM72}, \cite{Ke75} and \cite{So66}).  The preceding papers show the existence of {\em largest countable}  $\Sigma_n^1$ and $\Pi_n^1$ sets for certain $n$,  and also give information about their internal structure.  Incidentally they provide some nontrivial examples of non-typical reals. For instance,  Solovay \cite{So66} proved that if there are only countably many constructible reals (which is the case if there exists a measurable cardinal), then they form the largest countable $\Sigma_2^1$ set, which means that every constructible real is non-typical. Generalizing this, \cite{KM72} and  \cite{Ke75} showed that, under the assumption of Projective Determinacy (PD),  for each $n$ there exists  a largest countable $\Sigma_{2n}^1$ set, $C_{2n}$,  and a largest countable $\Pi_{2n+1}^1$ set, $C_{2n+1}$.

The above results   should perhaps prompt us to refine the general definition of typicality for reals, by considering   ``degrees of typicality'', namely $\Sigma_n$- and $\Pi_n$-typical reals, for $n\geq 0$, when the properties involved in the definitions are $\Sigma_n$ and  $\Pi_n$, respectively. The aim should be to find
conditions in order for a real to  belong to some countable $\Sigma_n$/$\Pi_n$ analytical set,  and  besides to specify ``concrete'' reals not belonging to any such set.

\begin{Rem} \label{R:nosets}
{\em Although our definitions of typical property and typical element apply to any kind of first-order structures, we shall refrain throughout from dealing with models of set theory, namely structures $\mbox{\goth M}=(M,E)$ for the language $L=\{\in\}$ that satisfy, say, the axioms of $ZF$. The reason is that  the definition of typical property over a structure $\mbox{\goth M}$ relies heavily on the} cardinality {\em of subsets of $M$ which is judged from} outside {\em the structure $\mbox{\goth M}$ (external cardinality). On the other hand, every model  $(M,E)$ of $ZF$ possess a  natural notion  of internal cardinality $|x|^{\mbox{\goth M}}$, for every $x\in M$, which provides a much more natural way to measure  the size of elements of $M$. Now if $\phi(x)$ is a property of the language $L=\{\in\}$, we cannot use the internal cardinality of  $\mbox{\goth M}=(M,E)$ to measure the extension of $\phi(x)$, since $\{x\in M:\mbox{\goth M}\models\phi(x)\}$ is a class in general, i.e., not an element  of $M$. Thus the criterion of majority extension cannot be directly applied to $\phi(x)$ when we employ internal cardinality. Nevertheless, we can get around this problem  if we use instead a parameterized version  of typicality,  $\alpha$-typicality, for every ordinal $\alpha\in M$, defined as follows:
A property $\phi(x)$ of the language of set theory (perhaps with parameters) is } $\alpha$-typical {\em over $\mbox{\goth M}$, if
$$\mbox{\goth M}\models |\{x\in V_{\alpha+1}:\phi(x)\}|>|\{x\in V_{\alpha+1}:\neg\phi(x)\}|.$$
Then $\phi(x)$ is said to be} typical over $\mbox{\goth M}$, {\em if it is $\alpha$-typical for every ordinal $\alpha\in M$.
This adjustment somehow internalizes the definition of typical property over $\mbox{\goth M}=(M,E)$, although  at the cost of increased complexity. But the internalization fails irreparably  when it comes to the definition of} $\alpha$-typical/typical element {\em  of $(M,E)$ (if in the old definition we simply replace ``typical property'' with ``$\alpha$-typical property''). For then  a set  $A\in M$ would be  $\alpha$-typical  iff  for every $\alpha$-typical property $\phi(x)$ over $\mbox{\goth M}$, $\mbox{\goth M}\models\phi(A)$.   The latter definition obviously cannot be expressed inside $\mbox{\goth M}$. Thus we are in a situation where one has to judge what a typical set of $\mbox{\goth M}$  is by using both internal and external criteria, a practice  not quite natural.

Despite the aforementioned subtleties and obstacles, the challenge to find  a  notion of typicality suitable for models of set theory remains. All the more because such models accommodate naturally} generic sets, {\em and it is well-known that  genericity is a notion dual to randomness in the contexts of both  reals and sets (see e.g. \cite{DH08}, \cite{La92}).
}
\end{Rem}

\section{Relative typicality as a relative randomness notion. Randomness axioms}

The notion of typicality considered above has a clear key-feature that allows it to be used as a faithful formal representation of the  (general) concept of randomness. This feature is the fact that  a typical object avoids all ``special'' predicates, i.e., predicates having  ``small extension'' (as measured by cardinality). All established   notions of randomness currently used  in the context of reals (namely ML-randomness, Schnorr randomness and  computable randomness), have similar features: a  random real  is one that avoids all sets of measure zero that are generated by certain special computable tests.

Because of the  discrepancy in the criteria of smallness and effectiveness, one should not expect to  prove many connections between typical reals and ML-random or Schnorr random reals. Rather the similarities should be sought in  the {\em structure} of the class of typical sets and that of the other classes of random sets. Such a structural approach to randomness, through axioms,  has been  set out  by van Lambalgen in a series of papers, \cite{La90},   \cite{La92} and \cite{La96}. In these papers the author introduced a new primitive  relation $R(x,\overline{y})$ for relative randomness, with  intended meaning ``$x$ is random with respect to $\overline{y}$'' (or ``$\overline{y}$ offers no information for $x$''), as well as  certain axioms   Ri about  $R$ (randomness axioms). In \cite{La90} the relation $R(x,\overline{y})$ refers primarily to  infinite binary sequences, i.e., elements of $2^{\omega}$, while in the sequel papers \cite{La92} and \cite{La96} $R(x,\overline{y})$ is intended to express properties of {\em general} random sets, so  the  axioms,  denoted again  Ri, are slightly different.

The relations $\textsf{Tp}(x,\overline{y})$ and $\textsf{wTp}(x,\overline{y})$ of relative typicality, defined in the previous section, present obvious analogies with  $R(x,\overline{y})$,  so it is natural to examine closely which of the axioms Ri hold of them and which do not.
Recall that when the relations $\textsf{Tp}$ and $\textsf{wTp}$ specialize to the language $L_2$ and the structure ${\cal R}$ of reals, they mean the following:

$\bullet$ $\textsf{Tp}(a,\overline{b}):\Leftrightarrow$ for every  formula $\phi(x,\overline{b})$ of $L_2$   such that   ${\cal R}\models\phi(a,\overline{b})$,   $|\{x:{\cal R}\models\phi(x,\overline{b})\}|=2^{\aleph_0}$.

$\bullet$ $\textsf{wTp}(a,\overline{b}):\Leftrightarrow$ for every  formula   $\phi(x,\overline{b})$  of $L_2$ such that   ${\cal R}\models\phi(a,\overline{b})$,   $|\{x:{\cal R}\models\phi(x,\overline{b})\}|\geq \aleph_0$.

The randomness axioms considered below are all  analogues of  axioms proposed  by van Lambalgen in \cite{La90}, \cite{La92}, and \cite{La96}. These papers  contain three  different  but largely overlapping lists of 6 or 7 axioms denoted Ri, for $i=1,\ldots,7$. The axioms  Ti, for $1\leq i\leq 6$, given below are the  analogues of the corresponding  Ri contained in  the list of   \cite{La90}. Axiom T7 (tailset) corresponds to R6 of \cite{La96}, while axiom-scheme T8 (zero-one law) corresponds to scheme R5 of  \cite{La96}.  Notice that the Ri's are formulated in a first-order object  language $L\cup\{R\}$, i.e., $L$ augmented by $R$. Moreover in the scheme  R6 and the zero-one law, the formulas may contain $R$. On the other hand  our Ti's are formulated in the metalanguage.   For finite tuples $\overline{a}$, $\overline{b}$, we write $\overline{a}\overline{b}$ for the concatenation of $\overline{a}$ and  $\overline{b}$.

It should be stressed that  not all axioms  Ti are expected to hold of the relation $\textsf{Tp}(x,\overline{y})$. The reason is that Ri were  mostly motivated by the intuition about  a random sequence of digits 0,1, as a sequence {\em generated by independent choices} (see \cite[p. 284]{La96}), while the intuition behind a typical object, as this is reflected in its definition,  is substantially different.

The next definition is needed for the formulation of axiom scheme T8 below.
\begin{Def} \label{D:tailset}
{\em The reals $a,b$ are said to be} equivalent {\em  and write $a\approx b$ if they differ in finitely many elements, i.e., if the symmetric difference $a\triangle b$ is finite. Let $[a]_\approx$ be the equivalence class of $a$. A set of reals $X$ is said to be a} tailset {\em  if $a\in X$ and $b\approx a$ imply $b\in X$.}
\end{Def}

\vskip 0.1in

T1. (Existence) $(\forall \overline{y})(\exists x)\textsf{Tp}(x,\overline{y})$. In particular, there exists $a$ such that $\textsf{Tp}(a)$.

T2. (Downward monotonicity) $\textsf{Tp}(a,c\overline{b})\Rightarrow \textsf{Tp}(a,\overline{b})$.

T3. (Parameters act as a  set rather than vector) (a) $\textsf{Tp}(a,\overline{b})\Rightarrow \textsf{Tp}(a,\pi\overline{b})$, for any permutation $\pi$ of $\overline{b}$. (b) $\textsf{Tp}(a,b\overline{c})\Rightarrow \textsf{Tp}(a,bb\overline{c})$.

T4. (Irreflexivity) $\textsf{Tp}(a,b)\Rightarrow a\neq b$.

T5. (Steinitz exchange principle) $\textsf{Tp}(a,\overline{c}) \ \& \ \textsf{Tp}(b,a\overline{c}) \Rightarrow \ \textsf{Tp}(a,b\overline{c})$.\footnote{A motivation for this principle comes from the  linear dependence relation.  As we said above the relation $R(a,\overline{b})$, and also $\textsf{Tp}(a,\overline{b})$,  can be read ``$a$ is independent from $\overline{b}$'', so $\neg \textsf{Tp}(a,\overline{b})$ can be seen as a dependence relation $D(a,\overline{b})$:  ``$a$ depends on  $\overline{b}$''. Then T5 is equivalently written  $D(a,b\overline{c})\Rightarrow D(a,\overline{c})\vee D(b,a\overline{c})$. This implication can be easily seen to be true if we take $a$, $b$, $c_1,\ldots,c_n$ to be vectors and $D(a,\overline{b})$ mean ``$a$ is a linear combination of $b_1,\ldots,b_n$.'' }

T6. If   $\phi(x,\overline{y})$ is a formula of $L$ not containing $z$ free, then
$$(\forall x)(\textsf{Tp}(x,z\overline{y})\Rightarrow \phi(x,\overline{y}))\Rightarrow (\forall x)(\textsf{Tp}(x,\overline{y})\Rightarrow \phi(x,\overline{y})).$$

T7. (Tailset, for the structure ${\cal R}$ of reals) For any reals  $\overline{c}$, the set $\{x:\textsf{Tp}(x,\overline{c})\}$ is a tailset.

T8. (Zero-one law, for the structure ${\cal R}$ of reals) Let $\phi(x,\overline{c})$ be a formula the extension of which is a tailset. Then

$(\exists x)(\textsf{Tp}(x,\overline{c})\ \& \ \phi(x,\overline{c}))$, implies

$(\forall x)(\textsf{Tp}(x,\overline{c})\Rightarrow  \phi(x,\overline{c}))$.

\begin{Prop} \label{P:randomaxiom}
Axioms {\rm T1}-{\rm T4} are true  in every structure $\mbox{\goth M}$. Axiom {\rm T7} holds in ${\cal R}$.
\end{Prop}

{\em Proof.}   T1: This follows immediately from  Theorem \ref{T:typicalreal}.
T2 and T3 are straightforward consequences of the definition of $\textsf{Tp}(a,\overline{b})$, since  what counts is simply the {\em set} of elements of  the tuple $\overline{b}$, not their particular ordering. T4 is equivalently written $\neg \textsf{Tp}(a,a)$, which follows immediately from Fact \ref{F:easy} (i).

T7: We have to show that for any $a,b,\overline{c}$, $\textsf{Tp}(a,\overline{c})$ and $a\approx b$ imply $\textsf{Tp}(b,\overline{c})$, or, equivalently, $\neg\textsf{Tp}(a,\overline{c})$ and $a\approx b$ imply $\neg\textsf{Tp}(b,\overline{c})$. Assume $\neg\textsf{Tp}(a,\overline{c})$ and $a\approx b$. Then there is $\phi(x,\overline{y})$ such that ${\cal R}\models\phi(a,\overline{c})$ and if $X=\{x:{\cal R}\models\phi(x,\overline{c})\}$, then $|X|<2^{\aleph_0}$. Notice that the predicate $Fin(x)$:= ``the real $x$ is finite'' is definable in ${\cal R}$, and thus so is the relation $x\approx y:=Fin(x\triangle y)$.  Put $\psi(x,\overline{c}):=(\exists y)(\phi(y,\overline{c})  \wedge x\approx y)$.  Then  ${\cal R}\models\psi(b,\overline{c})$.  It suffices to show that if $Y=\{x:{\cal R}\models\psi(x,\overline{c})\}$, then $|Y|<2^{\aleph_0}$. But   $Y=\bigcup\{[x]_\approx:{\cal R}\models \phi(x,\overline{c})\}$, and $|[x]_\approx|=\aleph_0$  for every real $x$. So $|Y|\leq |X|\cdot\aleph_0=|X|<2^{\aleph_0}$. \telos

\vskip 0.2in

Concerning the rest axioms T5, T6 and T8, they  either fail (in ${\cal  R}$ or in some other structure)  or their truth is open. We discuss separately the case of each one of them.

\vskip 0.1in

{\bf Axiom T5.} This axiom fails in  general  for both relations $\textsf{Tp}$ and $\textsf{wTp}$. It suffices to show in particular that
\begin{equation} \label{E:simple}
\textsf{Tp}(a) \ \& \ \textsf{Tp}(b,a) \ \Rightarrow \ \textsf{Tp}(a,b)
\end{equation}
is false over some structure $\mbox{\goth M}$, for some $a,b\in M$ (and $\overline{c}=\emptyset$). We show this  by the following counterexample.

\begin{Ex} \label{E:counter}
Let $\mbox{\goth M}=(M,S,+,\cdot,0)$ be a countable nonstandard model of $PA$ with the following property: There are $a,b\in M$ such that $a\notin {\rm df}(\emptyset)$, $a\in {\rm df}(b)$ and $b\notin {\rm df}(a)$. Then $\textsf{Tp}(a)$, $\textsf{Tp}(b,a)$ and $\neg \textsf{Tp}(a,b)$ hold in $\mbox{\goth M}$, so (\ref{E:simple}) is false over  $\mbox{\goth M}$.
\end{Ex}

{\em Proof.} First notice that there is an abundance of countable nonstandard models of $PA$ with the property stated above (e.g. every countable recursively saturated model has this property). Fix such an $\mbox{\goth M}$. Recall also (Fact \ref{F:easy})  that, because of  countability,  $\textsf{Tp}(a,b) \Leftrightarrow \textsf{wTp}(a,b)$ over $\mbox{\goth M}$, while (by Fact \ref{F:easy} again) $\textsf{wTp}(a,b)\Leftrightarrow a\notin {\rm cl}(b)$. In view of these equivalences it suffices  to show that there exist $a,b\in M$ such that
\begin{equation} \label{E:closure}
a\notin{\rm cl}(\emptyset) \ \& \ b\notin{\rm cl}(a) \ \& \ a\in {\rm cl}(b).
\end{equation}
By assumption there are $a,b\in M$ such that
$$a\notin{\rm df}(\emptyset) \ \& \ b\notin{\rm df}(a) \ \& \ a\in {\rm df}(b).$$
So in order to show (\ref{E:closure}) it suffices to establish that for every $\mbox{\goth M}\models PA$ and any $\overline{b}\in M$, ${\rm df}(\overline{b})={\rm cl}(\overline{b})$. But the latter is a well-known fact, analogous to Fact \ref{F:coincide},   that holds in every model of $PA$, due to the coding capabilities of integers.  \telos

\vskip 0.2in

However  the following remains open:

\begin{Que} \label{Q:T5}
Does {\rm T5}, and in particular  (\ref{E:simple}),  hold in the structure ${\cal R}$?
\end{Que}
Notice that if (\ref{E:simple}) were true over ${\cal R}$,  then any  for any typical reals  $a,b$, i.e., such that  $\textsf{Tp}(a)$ and  $\textsf{Tp}(b)$, we would have $\textsf{Tp}(a,b)\Leftrightarrow \textsf{Tp}(b,a)$. Some further open questions, analogous to \ref{Q:T5},  are cited at the end of the paper.

\vskip 0.1in

{\bf Axiom scheme T6.} The  explanation  for the corresponding scheme R6 (more precisely, for  its simpler  version R6$'$)  given in \cite{La90} is somewhat  obscure and hoc. (The formulation of and explanation for  R6 itself is much more complicated.)  In p.  1146 the author  justifies the introduction of R6$'$  as follows: for any $a,\overline{b}$ of a structure $\mbox{\goth M}$, if $a$ is random relative to  $\overline{b}$, i.e., if $R(a,\overline{b})$, then $a$ should not belong to the algebraic closure  of $\overline{b}$.  The author says  that this goal  is attained  when  the scheme R6$'$ is true, and this is the reason for adopting  it. He gives also a short proof that R6$'$ guarantees  the implication  $\mbox{\goth M}\models R(a,\overline{b}) \Rightarrow a\notin {\rm cl}(\overline{b})$.
In our case, however, where $R(a,\overline{b})$ is interpreted either as  $\textsf{Tp}(a,\overline{b})$ or as $\textsf{wTp}(a,\overline{b})$, it follows  {\em from the definition} (see Fact \ref{F:easy} (iv))  that
$\textsf{Tp}(a,\overline{b}) \ \Rightarrow \ \textsf{wTp}(a,\overline{b}) \ \Leftrightarrow \ a\notin {\rm cl}(\overline{b})$,
so no extra axiom is needed for this derivation. Apart from this,  it is very hard  to motivate  the truth of T6 for  the relations $\textsf{Tp}(a,\overline{b})$ and $\textsf{wTp}(a,\overline{b})$ over the reals. Besides, it is relatively easy to show that T6 fails over appropriate  countable non-standard models of $PA$. However  the proof is beyond the  main scope of the paper, and  we omit it. We conclude that  T6 is not a natural principle for typicality.

\vskip 0.1in

{\bf Axiom scheme T8.} Like T7, the scheme refers to reals again. Its truth for general formulas $\phi(x,\overline{y})$ is open to us. If for some formula $\phi(x,\overline{y})$ and some given parameters  $\overline{c}$, either  $\phi(x,\overline{c})$ or $\neg\phi(x,\overline{c})$ happens to be a {\em typical} property, then it is easy to see that T8 is true for such a  $\phi(x,\overline{c})$. Because  if $\phi(x,\overline{c})$ is typical, then {\em every} $\overline{c}$-typical element satisfies it, that is $(\forall x)(\textsf{Tp}(x,\overline{c})\Rightarrow  \phi(x,\overline{c}))$ is true.
If on the other hand $\neg\phi(x,\overline{c})$ is typical, then no $\overline{c}$-typical element satisfies $\phi(x,\overline{c})$, i.e., $(\exists x)(\textsf{Tp}(x,\overline{c})\ \& \ \phi(x,\overline{c}))$ is false. In both cases T8 holds true. However it is easy to find properties $\phi(x)$ with extension closed under $\approx$, such that neither $\phi(x)$ nor $\neg\phi(x)$ is typical. For instance $\phi(x):=(\forall n)(\exists m>n)(m \ \mbox{is even} \ \wedge \ m\in x)$ is such a property.

\subsection{Van Lambalgen's theorem}
In this section we compare typical reals mainly with ML-random reals, and also with Schnorr and computably random reals. For the precise definitions of these notions see \cite{DH08}.
In \cite{La87} van Lambalgen proved the following theorem which since then bears his name (see also \cite{DH08}, \S 11.6).

\begin{Thm} \label{T:Lambalgen}
For any  $a,b\subseteq \omega$, $a\oplus b$ is ML-random iff $a$ is ML-random and $b$ is $a$-ML-random.
\end{Thm}

Here $a\oplus b$ denotes  the set $\{2n:n\in a\}\cup\{2n+1:n\in b\}$. Does the relation $\textsf{Tp}(a,\overline{b})$  satisfy the analogue of Theorem \ref{T:Lambalgen} in the  structure ${\cal R}$? That is, is the equivalence $\textsf{Tp}(a\oplus b)\Leftrightarrow \textsf{Tp}(a) \  \& \ \textsf{Tp}(b,a)$ true for all reals $a,b$? We shall see that the answer is no. Instead we have the following simple reduction for the typicality of $a\oplus b$.

\begin{Thm} \label{T:oplus}
For any reals  $a$, $b$, $\overline{c}$,
$$\textsf{Tp}(a\oplus b,\overline{c})\Leftrightarrow \textsf{Tp}(a,\overline{c}) \vee \textsf{Tp}(b,\overline{c}).$$
\end{Thm}

{\em Proof.} We show equivalently
$$\neg\textsf{Tp}(a\oplus b,\overline{c})\Leftrightarrow \neg\textsf{Tp}(a,\overline{c}) \ \& \ \neg\textsf{Tp}(b,\overline{c}).$$

$\Rightarrow$: Assume $\neg\textsf{Tp}(a\oplus b,\overline{c})$. There is a $\phi(x,\overline{y})$ such that ${\cal R}\models\phi(a\oplus b,\overline{c})$ and $|\{x:{\cal R}\models\phi(x,\overline{c})\}|<2^{\aleph_0}$. Observe that, by the definition of the operation $\oplus$,  every real $x$ is  written in a unique way as a sum $x_0\oplus x_1$. So we can write $x_0=(x)_0$ and $x_1=(x)_1$, i.e. $x=(x)_0\oplus(x)_1$. Consider the formulas:

$\psi_0(x,\overline{c}):=(\exists y)(\phi(y,\overline{c}) \wedge x=(y)_0)$,

$\psi_1(x,\overline{c}):=(\exists y)(\phi(y,\overline{c}) \wedge x=(y)_1)$. \\
Since ${\cal R}\models\phi(a\oplus b,\overline{c})$, it follows that ${\cal R}\models\psi_0(a,\overline{c})$ and ${\cal R}\models\psi_1(b,\overline{c})$. Let $X=\{x:{\cal R}\models\phi(x,\overline{c})\}$, $X_0=\{x:{\cal R}\models\psi_0(x,\overline{c})\}$ and $X_1=\{x:{\cal R}\models\psi_1(x,\overline{c})\}$. It suffices to show that $|X_0|,|X_1|<2^{\aleph_0}$. Now using the functions $f_0(x)=(x)_0$ and $f_1(x)=(x)_1$, we have immediately that $f_0[X]=X_0$ and $f_1[X]=X_1$. Thus $|X_0|\leq |X|$ and $|X_1|\leq |X|$. Since by assumption $|X|<2^{\aleph_0}$, it follows that $|X_0|,|X_1|<2^{\aleph_0}$.

$\Leftarrow$: Assume $\neg\textsf{Tp}(a,\overline{c})$ and $\neg\textsf{Tp}(b,\overline{c})$. There are $\phi(x,\overline{c})$, $\psi(y,\overline{c})$ such that ${\cal R}\models\phi(a,\overline{c})$, ${\cal R}\models\psi(b,\overline{c})$, and if $X=\{x:{\cal R}\models\phi(x,\overline{c})\}$ and $Y=\{y:{\cal R}\models\psi(y,\overline{c})\}$, then $|X|, |Y|<2^{\aleph_0}$. Consider the formula
$$\sigma(z,\overline{c}):=(\exists x)(\exists y)(\phi(x,\overline{c}) \wedge \psi(y,\overline{c})\wedge z=x\oplus y).$$
Then ${\cal R}\models\phi(a,\overline{c})$ and  ${\cal R}\models\psi(b,\overline{c})$ imply  ${\cal R}\models\sigma(a\oplus b,\overline{c})$. Let $Z=\{z:{\cal R}\models\sigma(z,\overline{c})\}$. It suffices to show that $|Z|<2^{\aleph_0}$. Now since every $z\in Z$ is of the form $x\oplus y$ with $x\in X$ and $y\in Y$, it follows that $|Z|\leq |X\times Y|=|X|\cdot|Y|=\max(|X|,|Y|)<2^{\aleph_0}$. This completes the proof. \telos

\vskip 0.2in

Theorem \ref{T:oplus} makes the notion of typicality deviate considerably from ML-randomness. For example we have the following immediate consequence.

\begin{Cor} \label{C:deviate}
For every real $a$ and every non-typical $b$,
$\textsf{Tp}(a\oplus b)\Leftrightarrow \textsf{Tp}(a)$. In particular,
$\textsf{Tp}(a\oplus\emptyset)\Leftrightarrow \textsf{Tp}(a)$.
\end{Cor}

{\em Proof.} When  $\neg\textsf{Tp}(b)$, by \ref{T:oplus}, $\textsf{Tp}(a\oplus b)\Leftrightarrow \textsf{Tp}(a)\vee \textsf{Tp}(b)\Leftrightarrow \textsf{Tp}(a)$.  \telos

\begin{Cor} \label{C:fails}
Van Lambalgen's theorem \ref{T:Lambalgen} is false  for the relation $\textsf{Tp}(a,\overline{b})$ in ${\cal R}$. That is,  the equivalence $\textsf{Tp}(a\oplus b)\Leftrightarrow \textsf{Tp}(a) \  \& \ \textsf{Tp}(b,a)$ is false in general.
\end{Cor}

{\em Proof.} Otherwise we should have for every typical real $a$, $\textsf{Tp}(a\oplus \emptyset)\Leftrightarrow \textsf{Tp}(a) \  \& \ \textsf{Tp}(\emptyset,a)$, hence $\textsf{Tp}(\emptyset,a)$, which is obviously false. \telos

\vskip 0.2in

By \ref{C:deviate},  if $a$ is a typical set then so is  $a\oplus\emptyset=\{2n:n\in a\}$. This is in sharp contrast to ML-random sets. No set consisting of even (or odd) numbers alone, can be ML-random. This can be shown  easily using the definition of ML-random sets (see Theorem \ref{T:neg} below for the proof of a stronger fact),  but follows also immediately  from van Lambalgen's theorem \ref{T:Lambalgen}. For if $a\oplus\emptyset$ were ML-random then $a$ should be random and besides $\emptyset$ should be $a$-random, which is obviously false.

In my view the fact that a ML-random real  turns  into a non-ML-random one if we multiply its elements by 2, or by any constant $k$, is highly counterintuitive. Consequently so is also van Lambalgen's theorem which agrees with and is related to this fact. On the other hand, it is shown in  \cite{Yu07} that van Lambalgen's theorem fails for the two other randomness notions, namely Schnorr randomness and computable randomness. Specifically, direction  $\Rightarrow$ of Theorem \ref{T:Lambalgen} fails if ``ML-random'' is replaced by ``Schnorr random'' or ``computably random''.
Concerning the relationship between the three concepts, it is known that  if $\textbf{ML}$, $\textbf{SR}$, $\textbf{CR}$ denote the classes of these reals, respectively, then they are strictly nested as follows: $\textbf{ML}\varsubsetneq\textbf{CR}\varsubsetneq\textbf{SR}$.

Liang Yu (\cite{Yu07}) considers the  truth of van Lambalgen's theorem  so important, that he believes it should be the criterion as to which of the above mentioned  three randomness notions  should be finally accepted as the ``correct'' one. He concludes that this is Martin-L\"{o}f randomness since it is the only one that satisfies the theorem. He justifies his belief as follows:
\begin{quote}
``Philosophically, a random set should have the property that no information about any part of it can be obtained from another part. In particular, no information about `the left part' of a random set should be obtained from `the right part' and vice versa. In other words, `the left part' of a random set should be `the right part'-random and vice versa.''
\end{quote}
The problem with this argument is in   the meaning of the term ``part''. It is implicitly implied that any part of a random set should be random, and also random relative to any other part.  But why couldn't  a random set have  definable (or in general non-random) parts? I think that  a random set can have plenty of non-random parts.  According to my intuition about randomness, if we glue together a random set and a (disjoint) definable one, the outcome will be a random set. In general, if we glue together  a ``bad'' entity (random, irregular, undefinable, etc) with a ``good'' one (non-random, regular, definable), the ``bad'' entity prevails and spoils the composite whole.  (This general idea was examined  in  \cite{Tz15}, and was  shown that in many countable structures there are exist undefinable sets $X$ (called totally non-immune) with an abundance of definable parts, namely, for every definable set $A$, if $A\cap X$ is infinite, then it contains an infinite definable subset.)

Corollary \ref{C:deviate} helps answer the question about the comparability (with respect to inclusion) of the class of typical reals with the classes of ML-random, Schnorr random and computably random reals.  Let $\textbf{TP}$ be the class of typical reals.

\begin{Thm} \label{T:neg}
$\textbf{TP}$ is incomparable with  all classes $\textbf{ML}$, $\textbf{SR}$ and $\textbf{CR}$.
\end{Thm}

{\em Proof.} Since $\textbf{ML}\subseteq \textbf{CR}\subseteq \textbf{SR}$, it suffices to show that $\textbf{ML}\not\subseteq \textbf{TP}$ and $\textbf{TP}\not\subseteq \textbf{SR}$.  $\textbf{ML}\not\subseteq \textbf{TP}$  follows from the fact that there is a $\Delta_2^0$-definable ML-random real, namely Chaitin's $\Omega$ (see \cite{DH08}, \S 8.2). Then clearly $\Omega\in \textbf{ML}\backslash \textbf{TP}$.

To show $\textbf{TP}\not\subseteq \textbf{SR}$, pick a real  $a\in \textbf{TP}$. Then by \ref{C:deviate},  $a\oplus\emptyset\in \textbf{TP}$. It suffices to show that $a\oplus\emptyset\notin \textbf{SR}$. We saw already above that $a\oplus\emptyset\notin \textbf{ML}$, but this is not enough since  $\textbf{ML}\varsubsetneq\textbf{SR}$. We have to show that there is a Schnorr test $(V_n)_n$ such that $a\oplus\emptyset\in \bigcap_nV_n$. Recall that a Schnorr test is a ML-test  with the extra requirement that the measures of  $V_n$'s converge to zero in a {\em computable} way,  rather than a semi-computable one as is the case with the ML-tests. Namely, there must exist a computable strictly increasing $g:\omega\rightarrow \omega$ such that for every $n$, $\mu(V_n)=2^{-g(n)}$ (instead of simply $\mu(V_n)\leq 2^{-n}$). Since $a\oplus\emptyset=\{2n:n\in a\}$  consists of even numbers alone, we consider the  following sequence of sets. (In this argument we identify of course a set $a\subseteq \omega$ with its characteristic function.) For every $n$, let

$V_n=\{x\in 2^{\omega}:(\forall i\leq n)(x(2i+1)=0)\}$. \\
Then clearly  $a\oplus\emptyset\in \bigcap_nV_n$. So it remains to verify that $(V_n)_n$ satisfies the conditions required in order to be a Schnorr test.
Recall that for finite strings  $\sigma\in 2^{<\omega}$, the sets $[\sigma]=\{x\in 2^{\omega}:\sigma\subset x\}$ are the basic clopen sets of the topology of $2^{\omega}$. If $dom(\sigma)=\{0,\ldots,n-1\}$, we identify $\sigma$ with  the string $\sigma(0)\sigma(1)\cdots\sigma(n-1)$, and let $|\sigma|=n$ be the length of $\sigma$. Then $\mu([\sigma])=2^{-n}$.  It is easy to check that for every $n$,  the set $V_n$ above is written
$$V_n=\bigcup\{[i_00i_10i_2\cdots 0i_n0]:i_0,\ldots,i_n\in\{0,1\}\}.$$
Obviously $(V_n)_n$ is a computable sequence of open sets, each generated by finitely many basic sets. Namely, each  $V_n$ is the  union of the $2^{n+1}$ (disjoint) basic sets $[\sigma_m]=[i_00i_10i_2\cdots 0i_n0]$, where for every $m<2^{n+1}$,  $|\sigma_m|=2n+2$, thus $\mu([\sigma_m])=2^{-(2n+2)}$. Therefore
$$\mu(V_n)=\Sigma_{m< 2^{n+1}}\mu([\sigma_m])=2^{n+1}\cdot 2^{-(2n+2)}=2^{-(n+1)}.$$
Thus the sequence $(\mu(V_n))_n$ converges to zero in a  computable way, with $g(n)=n+1$. So $(V_n)_n$ is a Schnorr test. \telos

\subsection{Turing reducibility and jump operation}

The next theorem shows that the relation $\textsf{Tp}(a,\overline{c})$ is closed with respect to  Turing degrees, as well as with respect to the jump operation $a\mapsto a'$ (both ways).

\begin{Thm} \label{T:Turing}
For all reals $a$, $b$, $\overline{c}$ the following hold.

(i) $\neg\textsf{Tp}(a,\overline{c}) \ \& \ b\leq_T a \ \Rightarrow \ \neg\textsf{Tp}(b,\overline{c})$. (Equivalently: \\
$\textsf{Tp}(b,\overline{c}) \ \& \ b\leq_T a \ \Rightarrow \ \textsf{Tp}(a,\overline{c})$.)

(ii) $\textsf{Tp}(a,\overline{c}) \ \& \ b\equiv_T a \ \Rightarrow \ \textsf{Tp}(b,\overline{c})$.

(iii) $\textsf{Tp}(a,\overline{c}) \ \Leftrightarrow \ \textsf{Tp}(a',\overline{c})$, \\
where  $a'$ is the halting problem relative to $a$, i.e. $a'=\{e\in \omega:\Phi_e^a(e)\downarrow\}$.
\end{Thm}

{\em Proof.} (ii)  follows immediately from (i).

(iii) Since for every $a$, $a<_Ta'$, by (i),  $\neg\textsf{Tp}(a',\overline{c})\ \Rightarrow \neg\textsf{Tp}(a,\overline{c})$. For the converse, assume $\neg\textsf{Tp}(a,\overline{c})$. We shall show $\neg\textsf{Tp}(a',\overline{c})$. Then there is $\phi(x,\overline{c})$ such that ${\cal R}\models\phi(a,\overline{c})$ and if $X=\{x:{\cal R}\models\phi(x,\overline{c})\}$, then $|X|<2^{\aleph_0}$. Now since for every real $x$, $x'=\{e:\Phi^x_e(e)\!\downarrow\}$, it is rather straightforward  that there is a definable in ${\cal R}$ (without parameters) function $f$ on ${\cal P}(\omega)$ such that $f(x)=x'$. Consider the formula $\psi(x,\overline{c}):=(\exists y)(\phi(y,\overline{c}) \wedge x=f(y))$. Since ${\cal R}\models\phi(a,\overline{c})$ and $a'=f(a)$, it follows that  ${\cal R}\models\psi(a',\overline{c})$. Moreover, if $Y=\{x:{\cal R}\models\psi(x,\overline{c})\}$, clearly $Y=f[X]$. Thus $|Y|<2^{\aleph_0}$ since $|X|<2^{\aleph_0}$. Therefore $\neg\textsf{Tp}(a',\overline{c})$.

To prove (i) assume $\neg\textsf{Tp}(a,\overline{c})$ and $b\leq_T a$. Then there is a $\phi(x,\overline{y})$ such that ${\cal R}\models\phi(a,\overline{c})$ and $|\{x:\phi(x,\overline{c})\}|<2^{\aleph_0}$. Now the relation $x\leq_T y$ says: ``$x$ is  definable by a  $\Sigma_1^0$ formula that has $y$ as parameter'', so is definable in ${\cal R}$ by the help of  a $\Sigma_2^0$ satisfaction class for  all $\Sigma_1^0$ predicates. Let $\psi(x,y)$ denote the formula formalizing $x\leq_Ty$. Then  ${\cal R}\models \psi(b,a)$. Consider the formula
$$\sigma(x,\overline{c}):=(\exists y)(\phi(y,\overline{c}) \wedge \psi(x,y)).$$
Since ${\cal R}\models \phi(a,\overline{c})$ and $\psi(b,a)$, we have ${\cal R}\models\sigma(b,\overline{c})$. So in order to show that $\neg\textsf{Tp}(b,\overline{c})$, it suffices to show that $|\{x:{\cal R}\models\sigma(x,\overline{c})\}|<2^{\aleph_0}$.
If $X=\{y:{\cal R}\models \phi(y,\overline{c})\}$,  by assumption $|X|<2^{\aleph_0}$. Also for every real  $y$ the set $Y_y=\{x:\psi(x,y)\}=\{x:x\leq_Ty\}$ is countable. Therefore
$$|\{x:{\cal R}\models\sigma(x,\overline{c})\}|\leq \Sigma_{y\in X}|Y_y|=|X|\cdot\aleph_0=|X|<2^{\aleph_0}.$$ \telos

\subsection{Lower cones and some open questions}
In this last part we discuss some open questions about relative typicality of reals. Recall that the relation $\neg\textsf{Tp}(x,y)$ means something like ``$x$  depends on $y$'', but since  transitivity is (most probably) missing, this is  not a pre-ordering. On the other hand we saw in Theorem \ref{T:typicalreal} that for every $a$, the set $\{x:\neg\textsf{Tp}(x,a)\}$ is ``small'' (in fact a minority set, or just  countable under $CH$). This reinforces the intuition that the elements of   $\{x:\neg\textsf{Tp}(x,a)\}$, except of being relatively few, should not be   ``more complex'' than  $a$. Let us call the set $\{x:\neg\textsf{Tp}(x,a)\}$ {\em (lower) cone} of $a$, and denote it by $con(a)$, i.e.,
$$con(a)=\{x:\neg\textsf{Tp}(x,a)\}.$$
We  denoted by $\textbf{TP}$ the class of all typical reals. So ${\cal P}(\omega)\backslash \textbf{TP}=\{x:\neg \textsf{Tp}(x)\}$ is the class of non-typical ones. The following fact relates ${\cal P}(\omega)\backslash \textbf{TP}$ with the cones of reals.

\begin{Fac}  \label{F:nontyp}
${\cal P}(\omega)\backslash \textbf{TP}=con(\emptyset)= \bigcap\{con(a):a\in {\cal P}(\omega)\}$.
\end{Fac}

{\em Proof.} By axiom T2, for all $x,y$, $\neg\textsf{Tp}(x)\Rightarrow\neg\textsf{Tp}(x,y)$, so   ${\cal P}(\omega)\backslash \textbf{TP}\subseteq con(a)$ for every $a$. Therefore
$${\cal P}(\omega)\backslash \textbf{TP}\subseteq \bigcap\{con(a):a\in {\cal P}(\omega)\}\subseteq con(\emptyset).$$
On the other hand, by definition $\neg\textsf{Tp}(x)$ means $\neg\textsf{Tp}(x,\emptyset)$, so  ${\cal P}(\omega)\backslash \textbf{TP}=con(\emptyset)$. Thus the  three sets above coincide. \telos

\vskip 0.2in

The  intuition  expressed above on the relationship between the complexity of $a$ and that of the elements of $con(a)$, leads to the following question:
\begin{Que}  \label{Q:cone}
Is it true that $\neg\textsf{Tp}(a) \ \& \ b\in con(a) \Rightarrow \neg\textsf{Tp}(b)$?
\end{Que}

We have seen (Proposition \ref{P:randomaxiom}) that, as a consequence of axiom T7, $\textsf{Tp}(a) \ \& \ b\approx a \ \Rightarrow \ \textsf{Tp}(b)$. The following question seems also natural:

\begin{Que}  \label{Q:conetyp}
Does there exist a real $a$ such that  $\textsf{Tp}(a)$ but
$$(\forall x)(x\in con(a) \  \& \ b\not\approx a \ \Rightarrow \neg\textsf{Tp}(x))?$$
\end{Que}

Another question,  that we  already met before section 3.1 and was related to Question \ref{Q:T5},  can be formulated in terms of cones.

\begin{Que}  \label{Q:reform}
If $\textsf{Tp}(a)$ and $\textsf{Tp}(b)$, is it true that $b\in con(a)\Leftrightarrow a\in con(b)$?
\end{Que}

\vskip 0.2in

{\bf Acknowledgement} I would like to thank the anonymous referee for some helpful  remarks and suggestions.

\end{document}